\documentclass[12pt]{amsart}

\usepackage{amsthm, amsmath, amssymb, mathrsfs}
\usepackage{enumerate, microtype}
\usepackage{cite}

\usepackage[colorlinks]{hyperref}

\makeatletter
\@namedef{subjclassname@2010}{%
	\textup{2010} Mathematics Subject Classification}
\makeatother

\allowdisplaybreaks

\theoremstyle{plain}
\newtheorem{thm}{Theorem}[section]

\newtheorem{lem}[thm]{Lemma}
\newtheorem{prop}[thm]{Proposition}
\theoremstyle{definition}
\newtheorem{defn}[thm]{Definition}
\theoremstyle{remark}
\newtheorem{rem}[thm]{Remark}
\numberwithin{equation}{section}
\newtheorem*{nota}{Notation}

\title[Isomorphism Theorems for $PF_\Phi(G)$ and $PM_\Phi(G)$]{Isomorphism Theorems for the Algebras of $\Phi-$Pseudofunctions and $\Phi-$Pseudomeasures}
\author[A. Dabra]{Arvish Dabra$^\star$}
\address{Arvish Dabra,\newline\indent Department of Mathematics,\newline\indent Indian Institute of Technology Delhi,\newline\indent New Delhi - 110016, India.}
\email{arvishdabra3@gmail.com}
\author[N. S. Kumar]{N. Shravan Kumar}
\address{N. Shravan Kumar,\newline\indent Department of Mathematics,\newline\indent Indian Institute of Technology Delhi,\newline\indent New Delhi - 110016, India.}
\email{shravankumar.nageswaran@gmail.com}

\begin{document}


	\begin{abstract}
        In this article, we study the isomorphism problem for the algebras of $\Phi-$Pseudofunctions and $\Phi-$Pseudomeasures, denoted by $PF_\Phi(G)$ and $PM_\Phi(G),$ respectively. More precisely, for a certain class of Young functions $\Phi,$ we prove that if there exists an isometric isomorphism between $PF_\Phi(G_1)$ and $PF_\Phi(G_2),$ or between $PM_\Phi(G_1)$ and $PM_\Phi(G_2),$ then $G_1$ and $G_2$ are isomorphic as topological groups. In addition, we present an Orlicz version of Parrott's theorem.
	\end{abstract}
    
    
	\keywords{Banach algebras; Convolution operators; Isomorphism problem; Locally compact groups; Orlicz spaces.\\
		$\star$ Corresponding author. \\
		{\it Email address:} arvishdabra3@gmail.com (A. Dabra).}
	
	
	\maketitle

	
	\section{Introduction} 

    The isomorphism problem in Abstract Harmonic Analysis pertains to identify those groups that can recovered from one of their associated convolution algebras. The most extensively studied convolution algebras are the group algebra $L^1(G)$ and the measure algebra $M(G),$ where $G$ is any locally compact group.
    
    A foundational result in this direction is due to Wendel \cite{wendel}, who established that the Banach algebras $L^1(G_1)$ and $L^1(G_2)$ are isometrically isomorphic if and only if the underlying groups $G_1$ and $G_2$ are topologically isomorphic. In 1964, Johnson \cite{john} extended this line of research and proved an analogous result for the measure algebra $M(G).$ In 1968, motivated by Wendel's celebrated result, Parrott \cite[Theorem 2]{Par} proved that if there exists a surjective linear isometry $T:L^p(G_1) \to L^p(G_2),$ for $p \in (1,\infty)\setminus\{2\},$ such that $T(f \ast g) = T(f) \ast T(g)$ whenever $f,g,f \ast g \in L^p(G_1)$ and similarly for $T^{-1},$ then $G_1$ is isomorphic to $G_2$ as topological groups.

    A number of operator algebras are derived from the group algebra $L^1(G),$ most notably the reduced group $C^\ast$-algebra $C^\ast_r(G)$ and the group von Neumann algebra $VN(G).$ In the case of $VN(G),$ there exists examples of non-isomorphic groups whose associated von Neumann algebras are isometrically isomorphic. An example in this context is provided by the groups $G_1 = \mathbb{Z}_2 \oplus \mathbb{Z}_2$ and $G_2 = \mathbb{Z}_4,$ for which $VN(G_1) \cong VN(G_2)$ as von Neumann algebras, even though $G_1 \ncong G_2.$ 

    In 2022, Gardella and Thiel \cite{GH2} investigated the $L^p$-version of the isomorphism problem and extended the classical results of Wendel and Johnson to a broader class of convolution algebras, namely, $PF_p(G), PM_p(G)$ and $CV_p(G),$ which were introduced by Herz \cite{herz} several decades earlier. For $p \in [1,\infty) \setminus \{2\},$ they proved that if there exists an isometric isomorphism $CV_p(G_1) \cong CV_p(G_2)$ (or $PM_p(G_1) \cong PM_p(G_2)$ or $PF_p(G_1) \cong PF_p(G_2)$), then $G_1$ and $G_2$ are isomorphic as topological groups.

    It is well known that the Orlicz space $L^\Phi(G)$ is a significant generalization of the classical Lebesgue $(L^p)$ spaces, where $\Phi$ is a Young function. Over the past decade, Orlicz spaces have attracted considerable attention from various researchers. In 2023, Tabatabaie and Latifpour \cite{TL} studied the Orlicz $(L^\Phi)$-version of the isomorphism problem. They proved that, for a certain class of Young functions $\Phi,$ if $CV_\Phi(G_1)$ and $CV_\Phi(G_2)$ are isometrically isomorphic, then $G_1$ is topologically isomorphic to $G_2.$

    In this article, we investigate the isomorphism problem for the algebras $PF_\Phi(G)$ and $PM_\Phi(G).$ The structure of the paper is as follows. In Section \ref{sec2}, we recall the necessary definitions and notations. In Section \ref{sec3} and Section \ref{sec4}, we establish isomorphism theorems for the algebras $PF_\Phi(G)$ (Theorem \ref{mainthm}) and $PM_\Phi(G)$ (Theorem \ref{isomain}),  respectively. Specifically, for a certain class of Young functions $\Phi,$ we prove that if there exists an isometric isomorphism between $PF_\Phi(G_1)$ and $PF_\Phi(G_2),$ or between $PM_\Phi(G_1)$ and $PM_\Phi(G_2),$ then the underlying groups $G_1$ and $G_2$ are topologically isomorphic. Furthermore, in Section \ref{sec5}, we establish an Orlicz version of Parrott's theorem (Theorem \ref{parrott}).


    
    \section{Preliminaries}\label{sec2}

    We begin by recalling the basic terminology related to Orlicz spaces.

    An even convex function $ \Phi: \mathbb{R}\rightarrow [0,\infty]$ is called a Young function if it satisfies the conditions $\Phi(0)= 0$ and $\lim\limits_{x \to \infty} \Phi(x)= + \infty.$ Two Young functions $\Phi_1$ and $\Phi_2$ are said to be equivalent if there exists positive constants $k_1,k_2$ and $x_0 \geq 0$ such that
    $$\Phi_1(k_1 x) \leq \Phi_2(x) \leq \Phi_1(k_2 x), \hspace{1cm} \forall \, x \geq x_0.$$
     
    Given a Young function $\Phi,$ we define the convex function $\Psi$ by
    $$\Psi(y):= \sup{\{x\, |y|-\Phi(x):x\geq 0\}}, \hspace{1cm} y\in\mathbb{R}.$$
    This function $\Psi$ is also a Young function and is known as the complementary function to $\Phi.$ The pair $(\Phi,\Psi)$ is called a complementary pair of Young functions. As a standard example, for $1 < p < \infty,$ the function $\Phi(x) = |x|^p/p$ is a Young function whose complementary function is $\Psi(y) = |y|^q/q,$ where $1/p + 1/q = 1.$

    A continuous Young function $\Phi$ is called an $N$-function if it satisfies the following conditions: $\Phi(x) = 0$ if and only if $x=0,$ $\lim\limits_{x \to 0} \frac{\Phi(x)}{x} = 0$ and $\lim\limits_{x \to \infty} \frac{\Phi(x)}{x} = \infty.$ Throughout, we assume that $(\Phi,\Psi)$ is a complementary pair of $N$-functions.

    A Young function $\Phi$ is said to satisfy the $\Delta_{2}-$condition, denoted $\Phi \in \Delta_{2},$ if there exists a constant $k > 0$ and $x_{0} \geq 0$ such that $$\Phi(2x) \leq k \, \Phi(x),\hspace{1cm} \forall \, x \geq x_0.$$ A complementary pair of Young functions $(\Phi,\Psi)$ is said to satisfy the $\Delta_2-$condition if both $\Phi$ and $\Psi$ satisfy it individually. A Young function $\Phi$ is said to satisfy the $\Delta'-$condition, denoted $\Phi \in \Delta',$ if there exists a constant $k > 0$ and $x_{0} \geq 0$ such that 
    $$\Phi(x y) \leq k \, \Phi(x) \, \Phi(y), \hspace{1cm} \forall \, x,y \geq x_0.$$ It is easy to verify that the family of Young functions given by $\Phi_{\alpha,p}(x) = \alpha \, |x|^p,$ for $p \geq 1$ and $\alpha > 0,$ satisfies both the $\Delta_2-$condition and the $\Delta'-$condition.
    
    A Young function $\Phi$ is said to satisfy the Milnes-Akimoni$\check{c}$ condition (or MA condition) if, for each $\epsilon > 0,$ there exists constant $k_\epsilon > 1$ and $x_0(\epsilon) \geq 0$ such that $$\Phi'((1+\epsilon) \, x) \geq k_\epsilon \, \Phi'(x),\hspace{1cm} \forall \, x \geq x_0(\epsilon).$$ The family of $N$-functions given by $\Phi_\alpha(x) := |x|^\alpha (1 + |\log|x||),$ with $\alpha > 1,$ satisfies the MA condition. Moreover, the corresponding complementary pair of $N$-functions $(\Phi_\alpha,\Psi_\alpha)$ also satisfies the $\Delta_2-$condition.
    
    Given a Young function $\Phi,$ the associated Orlicz space is denoted by $L^{\Phi}(G)$ and defined as
    $${L}^{\Phi}(G) := \left\{ f: G \rightarrow  \mathbb{C}:f \, \mbox{is measurable and}\int_G\Phi(\alpha |f|)\ < \infty \text{ for some}~ \alpha>0  \right\}.$$ This space becomes a Banach space when equipped with the Luxemburg norm (also known as the Gauge norm), given by
    $$N_{\Phi}(f):= \inf \left\{k>0:\int_G\Phi\left(\frac{|f|}{k}\right) \leq1 \right\}.$$ Let $\Psi$ be the complementary function to $\Phi.$ The Orlicz norm $\|\cdot\|_{\Phi}$ on $L^\Phi(G)$ is then defined by $$\|f\|_{\Phi} := \sup \left\{\int_{G}|fg| :g\in L^\Psi(G) \, \, \text{and} \, \, \int_{G}\Psi(|g|) \leq1 \right\}.$$ These two norms are equivalent and for every $f \in L^\Phi(G),$ the following holds: $$N_\Phi(f) \leq \|f\|_\Phi \leq 2 \, N_\Phi(f).$$ 

    If the Young function $\Phi$ satisfies the $\Delta_{2}-$condition, then the space $\mathcal{C}_c(G)$ of all continuous functions on $G$ with compact support is dense in $L^\Phi(G).$ Moreover, if the complementary pair of Young functions $(\Phi,\Psi)$ satisfies the $\Delta_2-$condition, then $L^\Phi(G)$ is a reflexive Banach space \cite[Theorem 10, Pg. 112]{RR}.
    
    For further details on Orlicz spaces, we refer the readers to \cite{RR}.
    \newline

    Let $(\Phi,\Psi)$ be a complementary pair of Young functions satisfying the $\Delta_2-$condition. Let $M(G)$ denote the measure algebra on the locally compact group $G.$ For $\mu \in M(G)$ and $f \in L^\Phi(G),$ we define the left convolution operator $\lambda_\Phi(\mu): L^\Phi(G) \to L^\Phi(G)$ by $$(\lambda_\Phi(\mu))(g) := \mu \ast g, \hspace{1cm} (g \in L^\Phi(G)).$$ Let $\mathcal{B}(L^\Phi(G))$ be the Banach space of all bounded linear operators on $L^\Phi(G),$ equipped with the operator norm. It is straightforward to verify that $\lambda_\Phi(\mu) \in \mathcal{B}(L^\Phi(G))$ for every $\mu \in M(G).$ Define $PM_\Phi(G)$ as the closure of the set $\{\lambda_\Phi(\mu): \mu \in M(G)\}$ in $\mathcal{B}(L^\Phi(G))$ with respect to the ultra-weak topology (the $w^\ast$-topology). The algebra of $\Phi-$Pseudomeasures $PM_\Phi(G)$ is isometrically isomorphic to the dual space of the Orlicz Fig\`{a}-Talamanca Herz algebra $A_\Psi(G)$ \cite[Corollary 4.11]{AA}.
    
    If $\mu_f$ is the measure associated with $f \in L^1(G),$ then the norm closure of the set $\{\lambda_\Phi(f): f \in L^1(G)\}$ in $\mathcal{B}(L^\Phi(G))$ is denoted by $PF_\Phi(G)$ and termed as the algebra of $\Phi-$Pseudofunctions. It is easy to verify that $PM_\Phi(G)$ is the $w^\ast$-closure of $PF_\Phi(G)$ in $\mathcal{B}(L^\Phi(G)).$ Moreover, for each $\mu = \delta_a \, \, (\text{the Dirac measure at} \,\, a \in G),$ the corresponding operator $\lambda_\Phi(\delta_a)$ is denoted by $\lambda_\Phi(a).$ Thus, the map $\lambda_\Phi$ defines the left regular representation of $G$ on $L^\Phi(G).$

    Aghababa and Akbarbaglu \cite{AA} introduced the space of $\Phi-$convolution operators, denoted by $CV_\Phi(G),$ as the set of all bounded linear operators $T: L^\Phi(G) \to L^\Phi(G)$ that satisfy
    $$T(f \ast g) = T(f) \ast g,$$
    for all $f,g \in C_c(G).$ Note that $CV_\Phi(G)$ is a closed subalgebra of $\mathcal{B}(L^\Phi(G))$ and hence, a unital Banach algebra under composition. By \cite[Theorem 5.5]{AA}, an operator $T \in CV_\Phi(G)$ if and only if it commutes with right translations, i.e.,
    $$T(\rho_a(f)) = \rho_a(T(f)),$$
    for every $a \in G$ and $f \in L^\Phi(G),$ where $(\rho_a(f))(x) := f(xa)$ for $x \in G.$ Therefore, the definition of $CV_\Phi(G)$ by Aghababa and Akbarbaglu is equivalent to that of Tabatabaie and Latifpour \cite{TL}. It is well known that the following inclusions holds:
    $$PF_\Phi(G) \subseteq PM_\Phi(G) \subseteq CV_\Phi(G).$$

    For a comprehensive study of these Banach algebras, we refer the readers to the series of works \cite{AA, AD, ArRaSh, RLSK3, RLSK4}.
    \newline

    A locally compact group $G$ is said to be amenable if there exists a left-invariant mean on $L^\infty(G),$ i.e., a positive linear functional $\Lambda$ on $L^\infty(G)$ of norm one such that $$\langle \Lambda,L_a(f) \rangle = \langle \Lambda,f \rangle,$$ 
    for all $f \in L^\infty(G)$ and $a \in G,$ where $(L_a(f))(x) := f(a^{-1}x)$ denotes the left translation of $f.$ Classical examples of amenable groups include abelian groups, compact groups and solvable groups. However, the free group on two generators is not amenable.

    Throughout this article, $G$ is a locally compact group equipped with a fixed Haar measure and we assume that $(\Phi,\Psi)$ is a complementary pair of $N$-functions satisfying the $\Delta_2-$condition.



    \section{For the algebra of $\Phi-$Pseudofunctions $PF_\Phi(G)$}\label{sec3}
    
        In this section, we address the isomorphism problem for the Banach algebra $PF_\Phi(G).$ Our primary aim is to establish Theorem \ref{mainthm}, which shows that, for a certain class of Young functions $\Phi,$ the locally compact group $G$ can be recovered from the algebra $PF_\Phi(G).$ The ideas are inspired by Gardella and Thiel's work \cite{GH2} on $PF_p(G)$ and by the work of Tabatabaie and Latifpour \cite{TL} on $CV_\Phi(G).$ We start by introducing the relevant notations and definitions.

    \begin{nota}
        For any Banach algebra $\mathcal{A},$ we denote the left multiplier algebra of $\mathcal{A}$ by $\mathcal{M}_l(\mathcal{A}).$ The definition is as follows.
        $$\mathcal{M}_l(\mathcal{A}):= \{M \in \mathcal{B}(\mathcal{A}): M(a b) = M(a) \, b, \,\, \text{for all} \,\, a,b \in \mathcal{A}\}.$$
        From the definition, it is clear that $\mathcal{M}_l(\mathcal{A})$ is a closed subalgebra of $\mathcal{B}(\mathcal{A}).$
    \end{nota}

    Let us begin with the following proposition.

    \begin{prop}\label{inclusion}
        For any locally compact group $G,$ we have the following isometric inclusions:
        $$PF_\Phi(G) \subseteq \mathcal{M}_l(PF_\Phi(G)) \subseteq PM_\Phi(G) \subseteq CV_\Phi(G).$$
    \end{prop}

    \begin{proof}
        The first inclusion is evident from the fact that any Banach algebra $\mathcal{A}$ is \mbox{contained} in its left multiplier algebra $\mathcal{M}_l(\mathcal{A}).$ The final inclusion $PM_\Phi(G) \subseteq CV_\Phi(G)$ is well known. Further, by \cite[Theorem 4.1]{GH20}, there exists a unique isometric \mbox{representation} of $\mathcal{M}_l(PF_\Phi(G))$ as a unital subalgebra of $\mathcal{B}(L^\Phi(G)).$ With this \mbox{identification}, it is easy to verify that $\mathcal{M}_l(PF_\Phi(G))$ commutes with right \mbox{translations} and therefore lies in $CV_\Phi(G).$ 
        
        To prove the inclusion  $\mathcal{M}_l(PF_\Phi(G)) \subseteq PM_\Phi(G),$ let $M \in \mathcal{M}_l(PF_\Phi(G)).$ If $\{f_j\}_j$ is a bounded approximate identity for $L^1(G),$ then by \cite[Lemma 4.12]{TL}, the net $\{\lambda_\Phi(f_j)\}_j$ converges in the weak operator topology (WOT) to the identity operator $I_\Phi.$ Since the $w^\ast$-topology and the WOT coincide on bounded sets of $\mathcal{B}(E),$ where $E$ is a reflexive Banach space, it follows that $\{\lambda_\Phi(f_j)\}_j$ also converges to $I_\Phi$ in $w^\ast$-topology. As each $\lambda_\Phi(f_j) \in PF_\Phi(G),$ we have $M(\lambda_\Phi(f_j)) \in PF_\Phi(G)$ for all $j.$ Now, by using the fact that left multiplication on $CV_\Phi(G)$ is $w^\ast$-continuous, it follows that $\{M(\lambda_\Phi(f_j))\}_j$ converges to $M$  in $w^\ast$-topology. This implies
        $M \in PM_\Phi(G)$ and hence, the inclusion follows.
    \end{proof}

    \begin{nota}
        We denote by $M_\lambda^\Phi(G)$ the norm closure of $\lambda_\Phi(M(G))$ in $\mathcal{B}(L^\Phi(G)),$ where $M(G)$ is the measure algebra.
    \end{nota}

    \begin{rem}\label{multiplier}
        For $\mu \in M(G),$ we define the map $\lambda_\Phi(\mu): PF_\Phi(G) \to PF_\Phi(G)$ by
        $$(\lambda_\Phi(\mu))(\lambda_\Phi(f)):= \lambda_\Phi(\mu) \circ \lambda_\Phi(f) = \lambda_\Phi(\mu \ast f) \hspace{1cm} (f \in L^1(G)).$$
        It is straightforward to verify that $\lambda_\Phi(\mu) \in \mathcal{M}_l(PF_\Phi(G))$ and thus, we have $$M_\lambda^\Phi(G) \subseteq \mathcal{M}_l(PF_\Phi(G)).$$ Further, observe that this definition coincides with the identification of $\mathcal{M}_l(PF_\Phi(G))$ into $\mathcal{B}(L^\Phi(G)).$ 
    \end{rem}

    Recall that an element $v$ in a unital Banach algebra $\mathcal{A}$ is said to be an invertible isometry if it is invertible and $\|v\| = 1 = \|v^{-1}\|.$ 
    
    \begin{nota}
        We denote the group of invertible isometries in $\mathcal{A}$ by $\mathcal{U}(\mathcal{A}).$ Let $\mathcal{U}(\mathcal{A})_0$ be the connected component of $\mathcal{U}(\mathcal{A})$ in the norm topology containing the identity element of $\mathcal{A}.$ By \cite[Theorem 7.1, Pg. 60]{HR1}, $\mathcal{U}(\mathcal{A})_0$ is a closed normal subgroup of $\mathcal{U}(\mathcal{A}).$ Let $q$ be the canonical quotient map from $\mathcal{U}(\mathcal{A})$ onto $\mathcal{U}(\mathcal{A})/\mathcal{U}(\mathcal{A})_0.$ We denote the group $\mathcal{U}(\mathcal{A})/\mathcal{U}(\mathcal{A})_0$ by $q(\mathcal{U}(\mathcal{A})).$ 
    \end{nota}

    The following lemma determines the group of invertible isometries for the algebra $\mathcal{M}_l(PF_\Phi(G))$ for a certain class of Young functions $\Phi.$
    
    \begin{lem}\label{predual}
        Let $G$ be a locally compact group and $\Phi$ be a Young function which is not equivalent with the function $|\cdot|^2.$ Let $\Phi \in \Delta_2 \cap \Delta'$ and $L^\Phi(G) \subseteq L^1(G).$ Then the group of invertible isometries of $\mathcal{M}_l(PF_\Phi(G))$ coincides with that of $CV_\Phi(G).$
    \end{lem}

    \begin{proof}
        By Proposition \ref{inclusion}, it suffices to show that $\mathcal{U}(CV_\Phi(G)) \subseteq \mathcal{M}_l(PF_\Phi(G)).$ By \cite[Corollary 4.7]{TL}, this is further equivalent to proving that $\lambda_\Phi(a) = \lambda_\Phi(\delta_a)$ is in $\mathcal{M}_l(PF_\Phi(G))$ for all $a \in G,$ which is evident from Remark \ref{multiplier}.
    \end{proof}

    Recall that the strict topology (str) on $\mathcal{M}_l(PF_\Phi(G))$ is the restriction of the strong operator topology (SOT) on $\mathcal{B}(L^\Phi(G))$ to $\mathcal{M}_l(PF_\Phi(G)).$ 

    \begin{nota}
        If $\mathcal{A} = \mathcal{M}_l(PF_\Phi(G)),$ we write $\mathcal{U}(A)_{str}$ for the group $\mathcal{U}(A)$ equipped with the restriction of the strict topology. Similarly, we denote by $q(\mathcal{U}(A))_{str}$ for the group $q(\mathcal{U}(A))$ equipped with quotient topology induced by $\mathcal{U}(A)_{str} \to q(\mathcal{U}(\mathcal{A})).$ In fact, $W \subseteq q(\mathcal{U}(\mathcal{A}))$ is open if and only if $q^{-1}(W)$ is open in $\mathcal{U}(\mathcal{A})$ in the strict topology.
    \end{nota}

    The subsequent result is instrumental in establishing the main theorem of this section.

    \begin{thm}\label{isomor}
        Let $G$ be a locally compact group and $\Phi$ be a Young function which is not equivalent with the function $|\cdot|^2.$ Let $\Phi \in \Delta_2 \cap \Delta'$ and $L^\Phi(G) \subseteq L^1(G).$ Then the mappings
        $$\gamma_G : \mathbb{T} \times G \to \mathcal{U}(\mathcal{M}_l(PF_\Phi(G)))_{str} \,\,\, \text{and} \,\,\, \gamma'_G: G \to q(\mathcal{U}(\mathcal{M}_l(PF_\Phi(G))))_{str},$$
        given by $$\gamma_G(c,a) := c \, \lambda_\Phi(a) \,\,\, \text{and} \,\,\, \gamma'_G(a) := q(\lambda_\Phi(a)) = [\lambda_\Phi(a)],$$ are isomorphisms of topological groups.
    \end{thm}

    \begin{proof}
        By Lemma \ref{predual} and \cite[Theorem 4.5]{TL}, it follows that the maps $\gamma_G$ and $\gamma'_G$ are bijective group homomorphisms (ignoring the topologies). To establish that these maps are in fact topological isomorphisms, it suffices to prove that $\gamma_G$ is a homeomorphism. To show this, let $\{(c_j,a_j)\}_j$ be a net converging to $(c,a)$ in $\mathbb{T} \times G.$ Observe that, for $f \in L^1(G),$ 
        {\footnotesize{
        \begin{align*}
            \|(c_j \, \lambda_\Phi(a_j))(\lambda_\Phi(f)) - (c \, \lambda_\Phi(a))(\lambda_\Phi(f))\|
            &= \|c_j \, \lambda_\Phi(\delta_{a_j} \ast f) - c \, \lambda_\Phi(\delta_a \ast f)\|\\
            &\leq \|(c_j - c) \, \lambda_\Phi(\delta_{a_j} \ast f)\| + \|c \, (\lambda_\Phi(\delta_{a_j} \ast f) - \lambda_\Phi(\delta_a \ast f))\|\\
            &\leq |c_j-c| \, \|f\|_1 + |c| \, \|\delta_{a_j} \ast f - \delta_a \ast f\|_1.
        \end{align*}
        }}
        As $\delta_{a_j} \ast f \xrightarrow{\|\cdot\|_1} \delta_a \ast f$ in $L^1(G),$ it follows that
        $$(\gamma_G(c_j,a_j))(\lambda_\Phi(f)) \xrightarrow{\|\cdot\|_{PF_\Phi}} (\gamma_G(c,a))(\lambda_\Phi(f))$$
        in $PF_\Phi(G).$ Since $\lambda_\Phi(L^1(G))$ is norm dense in $PF_\Phi(G),$ we have
        $$\gamma_G(c_j,a_j) \xrightarrow{\text{str}} \gamma_G(c,a)$$ in $\mathcal{M}_l(PF_\Phi(G)).$

        For the converse, let $\{(c_j,a_j)\}_j$ be a net in $\mathbb{T} \times G$ and let $(c,a) \in \mathbb{T} \times G.$ Assume that $\gamma_G(c_j,a_j) \xrightarrow{str} \gamma_G(c,a)$ in $\mathcal{M}_l(PF_\Phi(G)),$ i.e.,
        $$(c_j \, \lambda_\Phi(a_j))(T) \xrightarrow{\|\cdot\|_{PF_\Phi}} (c \, \lambda_\Phi(a))(T)$$
        for all $T$ in $PF_\Phi(G).$

        Let $U$ be a neighborhood of the identity $e$ in $G.$ Then there exists a compact symmetric neighborhood $W$ of $e$ in $G$ such that $WWWW \subseteq U.$ Since $W$ is compact, we have $0 < |W| < \infty$ and thus, $\frac{\chi_W}{|W|} \in L^\Phi(G) \subseteq L^1(G).$ Observe that
        \begin{align*}
            c_j \, \chi_{a_jWW} &= \frac{c_j}{|W|} |a_jW| \, \chi_{a_jWW}
            = c_j \, \left(\chi_{a_jW} \ast \frac{\chi_W}{|W|}\right)\\
            &= (c_j \, \lambda_\Phi(\chi_{a_j W}))\left(\frac{\chi_W}{|W|}\right)
            = \left((c_j \, \lambda_\Phi(a_j))(\lambda_\Phi(\chi_W))\right)\left(\frac{\chi_W}{|W|}\right).
        \end{align*}
        Similarly,
        $$c \, \chi_{aWW} = \left((c \, \lambda_\Phi(a))(\lambda_\Phi(\chi_W))\right)\left(\frac{\chi_W}{|W|}\right).$$
        Thus, by taking $T = \lambda_\Phi(\chi_W) \in PF_\Phi(G)$ in the assumption, it follows that
        $$c_j \, \chi_{a_jWW} \xrightarrow{\|\cdot\|_\Phi} c \, \chi_{aWW}.$$
        Now, by applying the same argument as in \cite[Theorem 4.5]{TL} (with $W$ replaced by $WW$), it follows that $a_j$ eventually lies in $a \, U.$ Therefore, $a_j \to a$ in $G$ and hence, $c_j \to c$ in $\mathbb{T}.$ 
    \end{proof}

        We now present the main theorem of this section, which establishes that the Banach algebra $PF_\Phi(G)$ remembers $G.$ This result serves as an Orlicz analogue of \cite[Proposition 5.5]{GH2}.

    \begin{thm}\label{mainthm}
       Let $G_1$ and $G_2$ be two locally compact groups and $\Phi$ be a Young function which is not equivalent with the function $|\cdot|^2.$ Let $\Phi \in \Delta_2 \cap \Delta'$ and for $i =1,2,$ $L^\Phi(G_i) \subseteq L^1(G_i).$ If $PF_\Phi(G_1)$ and $PF_\Phi(G_2)$ are isometrically isomorphic as Banach algebras, then $G_1$ and $G_2$ are isomorphic as topological groups.
    \end{thm}

    \begin{proof}
        Let $\Gamma: PF_\Phi(G_1) \to PF_\Phi(G_2)$ be an isometric isomorphism. This induces an isometric isomorphism from $\mathcal{M}_l(PF_\Phi(G_1))$ to $\mathcal{M}_l(PF_\Phi(G_2)),$ denoted by $\tilde{\Gamma}$ and given by
        $$(\tilde{\Gamma}(M))(T) := \Gamma(M(\Gamma^{-1}(T))),$$
        for $M \in \mathcal{M}_l(PF_\Phi(G_1))$ and $T \in PF_\Phi(G_2).$ Since $\tilde{\Gamma}$ is an invertible isometry, it is a homeomorphism for both the norm topology and the strict topology.

        If we denote $\mathcal{M}_l(PF_\Phi(G_1))$ and $\mathcal{M}_l(PF_\Phi(G_2))$ by $\mathcal{A}_1$ and $\mathcal{A}_2,$ respectively, then $\tilde{\Gamma}$ induces a group isomorphism $\Gamma_\mathcal{U}:\mathcal{U}(\mathcal{A}_1) \to \mathcal{U}(\mathcal{A}_2).$ As $\Gamma_\mathcal{U}$ is a homeomorphism for the norm topology, it induces a group isomorphism $\Gamma_q: q(\mathcal{U}(\mathcal{A}_1)) \to q(\mathcal{U}(\mathcal{A}_2)).$ Further, as $\tilde{\Gamma}$ is also a homeomorphism for the strict topology, it follows that the map $\Gamma_\mathcal{U}: \mathcal{U}(\mathcal{A}_1)_{str} \to \mathcal{U}(\mathcal{A}_2)_{str}$ is an isomorphism of topological groups and hence, the group homomorphism $\Gamma_q: q(\mathcal{U}(\mathcal{A}_1))_{str} \to q(\mathcal{U}(\mathcal{A}_2))_{str}$ is also a topological isomorphism. Therefore,
        $$(\gamma'_{G_2})^{-1} \circ \Gamma_q \circ \gamma'_{G_1}:G_1 \to G_2$$
        is an isomorphism of topological groups, where $\gamma'_{G_1}: G_1 \to q(\mathcal{U}(\mathcal{A}_1))_{str}$ and $\gamma'_{G_2}: G_2 \to q(\mathcal{U}(\mathcal{A}_2))_{str}$ are topological group isomorphisms from Theorem \ref{isomor}. Hence, $G_1$ and $G_2$ are isomorphic as topological groups.
    \end{proof}

    \begin{rem}
        Recall that an isometric anti-isomorphism between two Banach \mbox{algebras} $\mathcal{A}_1$ and $\mathcal{A}_2$ is an isometric isomorphism from $\mathcal{A}_1$ to the opposite Banach algebra $\mathcal{A}_2^{op},$ where $\mathcal{A}_2^{op}$ shares the same underlying Banach space as $\mathcal{A}_2,$ but with the multiplication reversed. It is worth noting that Theorem \ref{mainthm} remains valid even when the Banach algebras $PF_\Phi(G_1)$ and $PF_\Phi(G_2)$ are isometrically anti-isomorphic.
    \end{rem}


    \section{For the algebra of $\Phi-$Pseudomeasures $PM_\Phi(G)$}\label{sec4}

    In this section, we investigate the isomorphism problem for the Banach algebra $PM_\Phi(G).$

    We begin by recalling the isomorphism theorem for the Banach algebra $CV_\Phi(G).$ Tabatabaie and Latifpour proved that for a certain class of Young functions $\Phi,$ if $CV_\Phi(G_1)$ is isometrically isomorphic to $CV_\Phi(G_2),$ then $G_1$ and $G_2$ are isomorphic as topological groups \cite[Theorem 4.13]{TL}.  

    One of the conditions satisfied by the aforementioned class of Young functions $\Phi$ is that $L^\Phi(G) \subseteq L^1(G).$ By \cite[Theorem 2]{hudzik}, this inclusion holds if and only if either the group $G$ is compact or the Young function $\Phi$ satisfies $$\lim_{x \to 0} \frac{\Phi(x)}{x} > 0.$$ However, we always assume that the pair $(\Phi,\Psi)$ is a complementary pair of $N$-functions. This implies that the Young function $\Phi$ satisfies $$\lim_{x \to 0} \frac{\Phi(x)}{x} = 0.$$ Therefore, the group $G$ must be compact.

    Since $G$ is compact and hence amenable, it follows from \cite[Theorem 6.4]{AA} that
    $$PM_\Phi(G) = CV_\Phi(G).$$ As a consequence of \cite[Theorem 4.13]{TL}, we obtain the following result.

    \begin{thm}\label{isomain}
        Let $G_1$ and $G_2$ be two locally compact groups and $\Phi$ be a Young function which is not equivalent with the function $|\cdot|^2.$ Let $\Phi \in \Delta_2 \cap \Delta'$ and for $i =1,2,$ $L^\Phi(G_i) \subseteq L^1(G_i).$ If $PM_\Phi(G_1)$ and $PM_\Phi(G_2)$ are isometrically isomorphic as Banach algebras, then $G_1$ and $G_2$ are isomorphic as topological groups.
    \end{thm}
    
    For $1 < p < \infty,$ Gardella and Thiel \cite[Theorem 6.3]{GH2} proved the reflexivity conjecture for all Banach subalgebras of $CV_p(G)$ that contain $PF_p(G).$ Motivated from this result, we prove the final theorem of this section, which characterizes the reflexive subalgebras of $CV_\Phi(G)$ containing $PF_\Phi(G)$ under the assumption that the group $G$ is amenable and $\Psi$ satisfies the MA condition.

    \begin{thm}
        Let $G$ be a locally compact group and let $(\Phi,\Psi)$ be a complementary pair of Young functions satisfying the $\Delta_2$-condition. Then any Banach algebra $\mathcal{A}$ such that $PF_\Phi(G) \subseteq \mathcal{A} \subseteq CV_\Phi(G)$ is finite dimensional if and only if $G$ is finite. Furthermore, if $G$ is amenable and $\Psi$ satisfies the MA condition, then $\mathcal{A}$ is reflexive if and only if $G$ is finite.
    \end{thm}

    \begin{proof}
        Assume that $G$ is finite. Then by \cite[Theorem 3.8]{AD}, $A_\Psi(G)$ is finite dimensional and consequently, its dual $PM_\Phi(G)$ is also finite dimensional. Since $G$ is finite, it is, in particular, amenable, compact and discrete. Therefore, by \cite[Corollary 4.13]{RLSK3} and \cite[Theorem 6.4]{AA}, we have
        $$PF_\Phi(G) = PM_\Phi(G) = CV_\Phi(G) = \mathcal{A}.$$
        Hence, the result follows.

        Conversely, suppose that $\mathcal{A}$ is finite dimensional. Then, being a closed subspace of $\mathcal{A},$ $PF_\Phi(G)$ is also finite dimensional. Since $A_\Psi(G)$ is a closed subspace of $W_\Psi(G) = PF_\Phi(G)^\ast$ \cite[Corollary 3.2 (iii)]{RLSK4}, it follows that $A_\Psi(G)$ is finite dimensional as well. Thus, by \cite[Theorem 3.8]{AD}, the group $G$ must be finite.

        Now, in the case of amenable groups, assume that $\mathcal{A}$ is reflexive. Then, by repeating the above arguments (replacing finite dimensional by reflexive), it follows that $A_\Psi(G)$ is reflexive. Therefore, by \cite[Corollary 3.9]{AD}, we conclude that $G$ is finite.
    \end{proof}


    \section{Orlicz version of Parrott's Theorem}\label{sec5}

    In the last section, the main objective is to prove Theorem \ref{parrott}, which can be regarded as an Orlicz analogue of Parrott's celebrated result \cite[Theorem 2]{Par}.
    
    \begin{defn}
        Let $f \in L^\Phi(G).$ We say that $f$ is a left multiplier if, for every $g \in L^\Phi(G),$ the convolution $f \ast g$ also belongs to $L^\Phi(G)$ and the corresponding left convolution operator, denoted by $\lambda_\Phi(f),$ is in $\mathcal{B}(L^\Phi(G)).$
    \end{defn}

        The following lemma serves as a key step in proving the main theorem of this section.

    \begin{lem}\label{weakstar}
        Let $G$ be a locally compact group and let $\Phi \in \Delta_2.$ If $\mathcal{A}$ denotes the $w^\ast$-closure of the set $\{\lambda_\Phi(f):f \in L^\Phi(G) \,\, \text{is a left multiplier}\}$ within $\mathcal{B}(L^\Phi(G)),$ then $\mathcal{A}$ is a $w^\ast$-closed subalgebra of $\mathcal{B}(L^\Phi(G))$ such that $PM_\Phi(G) \subseteq \mathcal{A} \subseteq CV_\Phi(G).$
    \end{lem}

    \begin{proof}
        Let us denote the set $\{\lambda_\Phi(f):f \in L^\Phi(G) \,\, \text{is a left multiplier}\}$ by $\mathcal{A}_0.$ Then, by the definition of a left multiplier, it is evident that $\mathcal{A}_0$ forms a subalgebra of $\mathcal{B}(L^\Phi(G)).$ Hence, by assumption, $\mathcal{A}$ is a $w^\ast$-closed subalgebra of $\mathcal{B}(L^\Phi(G)).$ 
        To show that $\mathcal{A} \subseteq CV_\Phi(G),$ let $a \in G.$ Then, for all $g \in L^\Phi(G),$ we have
        $$(\lambda_\Phi(f))(\rho_a(g)) = f \ast (\rho_a(g)) = \rho_a(f \ast g) = \rho_a((\lambda_\Phi(f))(g)).$$
        This implies that $\lambda_\Phi(f)$ commutes with right translations and thus, the inclusion is clear.
        Furthermore, observe that $\lambda_\Phi(C_c(G)) \subseteq \mathcal{A}_0,$ since each $f \in C_c(G)$ is in $L^\Phi(G)$ and $$\|(\lambda_\Phi(f))(g)\|_\Phi = \|f \ast g\|_\Phi \leq \|f\|_1 \|g\|_\Phi \hspace{1cm} (g \in L^\Phi(G)),$$
        which shows that the corresponding left convolution map is bounded. Consequently,
        $$\overline{\lambda_\Phi(C_c(G))}^{w^\ast} \subseteq \overline{\mathcal{A}_0}^{w^\ast} = \mathcal{A}.$$ 
        Since $C_c(G)$ is dense in $L^1(G)$ in the norm topology, the inclusion $PM_\Phi(G) \subseteq \mathcal{A}$ follows.
    \end{proof}

    \begin{rem}\label{remmulti}
        If $L^\Phi(G) \subseteq L^1(G),$ then by \cite[Theorem 3.1]{oztop1}, $L^\Phi(G)$ is a convolution algebra and hence, $\lambda_\Phi(f)$ is a left multiplier for all $f \in L^\Phi(G).$
    \end{rem}

        We are now ready to prove the Orlicz version of Parrott's theorem.

    \begin{thm}\label{parrott}
        Let $G_1$ and $G_2$ be two locally compact groups and $\Phi$ be a Young function which is not equivalent with the function $|\cdot|^2.$ Let $\Phi \in \Delta_2 \cap \Delta'$ and for $i =1,2,$ $L^\Phi(G_i) \subseteq L^1(G_i).$ If there exists a surjective linear isometry $T: L^\Phi(G_1) \to L^\Phi(G_2)$ such that $$T(f \ast g) = T(f) \ast T(g) \,\,\hspace{2.17cm} \text{for all} \,\,\, f,g \in L^\Phi(G_1)$$ and $$T^{-1}(f \ast g) = T^{-1}(f) \ast T^{-1}(g) \hspace{1cm} \text{for all} \,\,\, f,g \in L^\Phi(G_2),$$ then $G_1$ and $G_2$ are isomorphic as topological groups.
    \end{thm}

    \begin{proof}
        For $i =1,2,$ define $$\mathcal{A}_0^i := \{\lambda_\Phi(f):f \in L^\Phi(G_i)\}.$$ Let $\mathcal{A}_i$ denote the $w^\ast$-closure of $\mathcal{A}_0^i$ in $\mathcal{B}(L^\Phi(G_i))$ for each $i.$ Then, by Lemma \ref{weakstar} and Remark \ref{remmulti}, it follows that $\mathcal{A}_i$ is a $w^\ast$-closed subalgebra of $\mathcal{B}(L^\Phi(G_i))$ satisfying 
        $$PM_\Phi(G_i) \subseteq \mathcal{A}_i \subseteq CV_\Phi(G_i),$$
        for each $i=1,2.$ As noted in Section \ref{sec4}, each $G_i$ must be compact and therefore, we have
        $$\mathcal{A}_i = PM_\Phi(G_i) = CV_\Phi(G_i),$$ for $i =1,2.$ Let $\Gamma: \mathcal{B}(L^\Phi(G_1)) \to \mathcal{B}(L^\Phi(G_2))$ be given by
        $$\Gamma(S) := T \circ S \circ T^{-1},$$
        for $S \in \mathcal{B}(L^\Phi(G_1)).$ It is easy to verify that $\Gamma$ is an isometric isomorphism of Banach algebras. Let $\Gamma_r$ denote the restriction of $\Gamma$ to $\mathcal{A}_1$ and let $f \in L^\Phi(G_1).$ Then, for all $g \in L^\Phi(G_2),$
        \begin{align*}
        (\Gamma(\lambda_\Phi(f)))(g) &= (T \circ \lambda_\Phi(f) \circ T^{-1})(g) = (T \circ \lambda_\Phi(f))(T^{-1}(g))\\  &= T(f \ast T^{-1}(g)) = T(f) \ast g = (\lambda_\Phi(T(f)))(g).
        \end{align*}
        It follows from the previous computation that $\Gamma_r$ maps $\mathcal{A}_0^1$ into $\mathcal{A}_0^2.$ Since $\Gamma_r$ is also $w^\ast$-continuous, this extends to a map from $\mathcal{A}_1$ to $\mathcal{A}_2,$ i.e., $\Gamma_r : \mathcal{A}_1 \to \mathcal{A}_2.$ Applying the similar arguments to $T^{-1}$ shows that $\Gamma_r$ is an isometric isomorphism between $\mathcal{A}_1$ and $\mathcal{A}_2.$ Therefore, by Theorem \ref{isomain}, it follows that $G_1$ and $G_2$ are isomorphic as topological groups.
    \end{proof}
   
    
    \section*{Acknowledgement}
    The first author gratefully acknowledges the Indian Institute of Technology Delhi for providing the Institute Assistantship.

    \section*{Data Availability} 
    Data sharing does not apply to this article as no datasets were generated or analysed during the current study.

    \section*{Competing Interests}
    The authors declare that they have no competing interests.   
    
	\bibliographystyle{acm}
	\bibliography{ref}
	
\end{document}